\documentclass[12pt]{article}
\usepackage{amsmath,amsthm,amsfonts}

\newtheorem{thm}{Theorem}[section]
\newtheorem{lemma}[thm]{Lemma}

\newtheorem{cor}[thm]{Corollary}
\newtheorem{defin}[thm]{Definition}
\newtheorem{rem}[thm]{Remark}
\newtheorem{exam}[thm]{Example}

\newcommand{\R}{{\mathbb{R}}}

\newcommand{\Z}{{\mathbb{Z}}}

\newcommand{\C}{{\mathbb{C}}}

\newcommand{\cD}{{\mathcal{D}}}

\newcommand{\Ca}{{\hbox{\it Cal}}}

\newcommand{\Ham}{{\hbox{\it Ham\,}}}
\newcommand{\Symp}{{\hbox{\it Symp\,}}}

\newcommand{\id}{{\text{{\bf 1}}}}

\newcommand{\Qed}{\hfill \qedsymbol \medskip}

\begin{document}

\title{Calabi quasimorphisms for the symplectic ball}

\author{ \textsc Paul Biran, Michael Entov$^1$\ and Leonid
  Polterovich$^2$}

\date{29 October 2003}

\maketitle

\footnotetext[1]{Partially supported by the Fund for Promotion of
  Research at Technion and by the Israel Science Foundation grant $\#$
  68/02.}

\footnotetext[2]{Supported by the Israel  Science Foundation
founded by the Israel Academy of Sciences and Humanities}

\begin{abstract}
   We prove that the group of compactly supported symplectomorphisms
   of the standard symplectic ball admits a continuum of linearly
   independent real-valued homogeneous quasimorphisms.  In addition
   these quasimorphisms are Lipschitz in the Hofer metric and have the
   following property: the value of each such quasimorphism on any
   symplectomorphism supported in any "sufficiently small" open subset
   of the ball equals the Calabi invariant of the symplectomorphism.
   By a "sufficiently small" open subset we mean that it can be
   displaced from itself by a symplectomorphism of the ball. As a
   byproduct we show that the (Lagrangian) Clifford torus in the
   complex projective space cannot be displaced from itself by a
   Hamiltonian isotopy.
\end{abstract}

\vfill\eject


\section{Introduction and results} 

A quasimorphism on a group $G$ is a function $\mu: G \to \R$ which
satisfies the homomorphism equation up to a bounded error: there
exists $R > 0$ such that $$|\mu(fg) -\mu(f) -\mu(g)| \leq R$$
for all
$f,g \in G$ (see \cite{Bav} for preliminaries on quasimorphisms). A
quasimorphism $\mu$ is called {\it homogeneous} if $\mu (g^m) = m \mu
(g)$ for all $g \in G$ and $m \in \Z$. It is easy to see that a
homogeneous quasimorphism on an abelian group is always a genuine
homomorphism.

In this paper we focus on the case when $G = \Symp_0 (B^{2n})$ is the
identity component of the group of all compactly supported symplectic
diffeomorphisms of the $2n$-dimensional ball
$$B^{2n} = \{|p|^2 +|q|^2 < 1\} \subset \R^{2n}$$
equipped with the
symplectic form $dp \wedge dq$, where $p = (p_1,\ldots, p_n),
q=(q_1,\ldots, q_n)$ are coordinates on $\R^{2n}$.  A celebrated
result due to A.Banyaga \cite{Ban} states that the real vector space
of homomorphisms $G \to \R$ is one-dimensional: all of them are
proportional to the classical Calabi homomorphism (see Section
\ref{Ham-etc} below).  Our main result shows that in contrast to this
$G$ carries a lot of homogeneous quasimorphisms.

\begin{thm}\label{thm-continuum-quasim-ball-0}
   The real vector space of homogeneous quasimorphisms on
   $\Symp_0(B^{2n})$ is infi\-nite-dimensional.
\end{thm}

\begin{cor}\label{cor-bd-hom}
   The second bounded cohomology of $\Symp_0\, (B^{2n})$ is an
   infi\-nite-dimensional vector space over $\R$.
\end{cor}

It is unknown whether analogues of these results are valid for
smooth and/or volume-preserving diffeomorphisms. The property
featured in the theorem above is known for Gromov-hyperbolic
groups \cite{BF,EF}. We refer the reader to \cite{EP} for a
historical review of the "hunt" for quasimorphisms. Let us mention
that the first genuine homogeneous quasimorphism (i.e. not a
homomorphism) on $\Symp_0 (B^{2n})$ was discovered by Barge and
Ghys in \cite{Ba-Ghys}.

The case $n=2$ was settled by Entov-Polterovich in \cite{EP} and by
Gambaudo-Ghys in \cite{Gamb-Ghys}. Though the methods used in these
papers are quite different, both of them appeal in a crucial way to
2-dimensional topology. To prove Theorem
\ref{thm-continuum-quasim-ball-0} we extend the symplecto-topological
approach developed in \cite{EP}.

The quasimorphisms guaranteed by
Theorem~\ref{thm-continuum-quasim-ball-0} have a number of
interesting additional features which will be described below:
they are Li\-pschitz with respect to the Hofer metric \cite{Hofer,
Pol-book} and enjoy the so-called Calabi property.

As a byproduct we prove a result on Lagrangian intersection.  Consider
the complex projective space ${\mathbb{C}}P^n$ endowed with its
standard symplectic structure. Then the so called {\it Clifford torus}
\[
T^n_{clif} := \{ [z_0:\ldots : z_n]\in \C P^n \ | \ \ |z_0| =\ldots =
|z_n|\} \subset \C P^n,
\]
is a Lagrangian torus (see Section~\ref{qm-cliff}).  Note that
under the usual symplectic identification of ${\mathbb{C}}P^n
\setminus {\mathbb{C}}P^{n-1}$ with the open unit ball
$\textnormal{Int\,}B^{2n}(1)$, the torus $T^n_{clif}$ corresponds
to the "split" torus $\{ (w_1,\ldots, w_n)\in\C^n \, | \ \ |w_i|^2
= \frac{1}{n+1},\ i=1,\ldots, n\}$ (see
Section~\ref{subsection-pf-thm-continuum-quasim-ball}).

\begin{thm}
   \label{thm-cliff-torus-not-displaceable}
   The Lagrangian torus $T^n_{clif}\subset \C P^n$ cannot be displaced
   from itself by a Hamiltonian isotopy.
\end{thm}

This result admits a generalization to certain products of complex
projective spaces (see Example~\ref{ex-1} in
Section~\ref{sec-outl}).

\section{Hamiltonian diffeomorphisms and Calabi quasimorphisms}
\label{Ham-etc}
We will work in the following general setting.  Let $(M, \omega)$ be a
connected symplectic manifold without boundary and let $G =
\Ham(M,\omega)$ denote the group of compactly supported Hamiltonian
symplectomorphisms of $(M,\omega)$.  For simply connected symplectic
manifolds $M$ (in particular, balls) $G$ coincides with the identity
component of the group of all compactly supported symplectomorphisms
of $M$.  For preliminaries on $G$ see e.g.
\cite{McD-Sal-sympl-top},\cite{Pol-book}.

Given a Hamiltonian $H: [0,1]\times M\to\R$ set $H_t := H(t, \cdot ):
M\to \R$ and denote by $\{ \phi_H^t\}_{0\leq t \leq 1}$, $\phi_H^0
=\id$, the Hamiltonian flow of $H$. Denote the time-1 map of this flow
by $\phi_H\in G$.

The group $G$ has a natural class of subgroups associated to non-empty
open subsets $U \subseteq M$. Namely, consider Hamiltonian functions
$H: [0,1]\times M\to\R$ such that for each $t$ the support of $H_t$
lies inside $U$. The subgroup $G_U\subset G$ is formed by all elements
$\phi_H$ generated by such Hamiltonian functions $H$. When the
symplectic form $\omega$ is exact on $U$ the formula
\begin{equation}
   \label{def-calabi-homom} \phi_H \mapsto \int_0^1 dt \int_M H_t
   \omega^n.
\end{equation}
gives rise to a well-defined homomorphism
\[
\Ca_U: G_U\to\R,
\]
called the {\it Calabi homomorphism} \cite{Ban}, \cite{Cal}. Note that
$G_U \subset G_V$ for $U \subset V$ and in this case $\Ca_U = \Ca_V$
on $G_U$.

In what follows we deal with the class $\cD_{ex}$ of all non-empty
open subsets $U$ which can be {\it displaced} by a Hamiltonian
diffeomorphism:
\[
hU \cap \text{Closure}~(U) = \emptyset \;\text{for}\;\text{some}\; h
\in G,
\]
and such that $\omega$ is exact on $U$. The latter condition is
automatically fulfilled when $M$ is a subset of the standard
symplectic $\R^{2n}$.

\begin{defin}
   \label{def-Calabi-quasi} {\rm A quasimorphism on $G$ coinciding
     with the Calabi homomorphism $\{\Ca _U: G_U \to \R\}$ on every ${U
       \in \cD_{ex}}$ is called {\it a Calabi quasimorphism}. }
\end{defin}

\begin{defin}
\label{def-lip}{\rm We say that a quasimorphism
  $\mu$ is {\it Lipschitz with respect to the Hofer metric} if there
  exists a constant $K > 0$ so that
  $$|\mu(\phi_F) -\mu(\phi_H)| \leq K \cdot ||F-H||_{C^0}$$
  for all
  Hamiltonians $F,H:[0;1]\times M \to \mathbb{R}.$ (For the relation
  of $||F-H||_{C^0}$ to the Hofer distance between $\phi_F$ and
  $\phi_H$ see \cite{Hofer, Pol-book}, also see \cite{EP}).  }
\end{defin}

With this language we can refine Theorem
\ref{thm-continuum-quasim-ball-0} as follows.

\begin{thm}\label{thm-continuum-quasim-ball-1}
   The real affine space of homogeneous Calabi quasi\-morph\-isms on
   $\Symp_0(B^{2n})$ which are Lipschitz with respect to the Hofer
   metric is infinite-dimensional.
\end{thm}

\begin{rem}{\rm  A straightforward modification of our arguments
     below \break yields the same result when one replaces the ball $B^{2n}$
     by the unit ball bundle over the $n$-torus, $T^n \times B^n
     \subset T^*T^n$, equipped with the standard symplectic
     structure.}
\end{rem}

The key ingredient in our proof of Theorem
\ref{thm-continuum-quasim-ball-1} is the following fact established in
\cite{EP}:
\begin{thm} \label{thm-existence}
   The group of Hamiltonian diffeomorphisms of the complex projective
   space $\C P^n$ endowed with the Fubini-Study symplectic form
   carries a homogeneous Calabi quasimorphism which is Lipschitz with
   respect to the Hofer metric.
\end{thm}
The paper \cite{EP} contains an explicit construction of such a
quasimorphism. However it is unknown whether it is unique.

\section{A Calabi quasimorphism on ${\mathbb{C}}P^n$ and the Clifford torus}
\label{qm-cliff}

Consider $M=\C P^n$ equipped with the Fubini-Study symplectic
structure $\omega$ normalized so that the integral of $\omega$ over
the complex projective line is 1 and hence $\textnormal{Vol}\, (\C
P^n) := \int_{{\mathbb{C}}P^n} \omega^{ n} = 1$.  The torus $T^n =
\R^n / \Z^n$ acts on $\C P^n$ in a Hamiltonian way as follows:
\[
(e^{2\pi s_1},\ldots, e^{2\pi s_n}): \ \ [z_0 :\ldots : z_n] \mapsto
[z_0 : e^{2\pi s_1} z_1 :\ldots : e^{2\pi s_n} z_n].
\]
A moment map $\Phi: \C P^n\to\R^n$ for the action is given by the
formula:
\[
\Phi : \ \ [z_0 :\ldots : z_n] \mapsto {\Big (}
\frac{|z_1|^2}{|z_0|^2+\ldots+|z_n|^2}, \ldots ,
\frac{|z_n|^2}{|z_0|^2+\ldots+|z_n|^2} {\Big )}.
\]
The image of $\Phi$ is the following closed convex polytope:
\[
\Delta= \{ p\in \R^n\ | \ p_1 +\ldots+p_n \leq 1\, ; \, p_i\geq 0, \,
i=1,\ldots, n \},
\]
where $p=(p_1,\ldots, p_n)$ are coordinates on $\R^n$. The inverse
image of a point $p\in\Delta$ under $\Phi$ is an isotropic torus which
is an orbit of the $T^n{\hbox{\rm -action}}$ and whose dimension is
equal to $i$ if and only if $p$ lies inside some $i{\hbox{\rm
    -dimensional}}$ face of $\Delta$ (i.e. it belongs to that
$i{\hbox{\rm -dimensional}}$ face but not to any face of smaller
dimension). Let
\[
p_{clif} := {\Big (}\frac{1}{n+1},\ldots, \frac{1}{n+1}{\Big )}
\]
denote the barycenter of $\Delta$.  The Lagrangian torus $T^n_{clif}:=
\Phi^{-1} (p_{clif})$ is called the {\it Clifford torus}.

We say that a function on the simplex $\Delta$ is smooth if it extends
to a smooth function in a neighborhood of $\Delta$ in $\mathbb{R}^n$.
The proof of our main results is based on the following calculation.

\begin{thm}
   \label{thm-calabi-clifford-cpn}
   Let $\mu: G \to \mathbb{R}$ be any Calabi quasimorphism which is
   continuous with respect to the Hofer metric. Then for every
   autonomous Hamiltonian $F:\C P^n\to \R$ of the form $\Phi^\ast
   \bar{F}$, where $\bar{F}: \Delta\to\R$ is a smooth function on
   $\Delta$,
   \[
   \mu (\phi_F) = \int_{\C P^n} F\omega^n - \bar{F} (p_{clif}).
   \]
\end{thm}

The proof of Theorem~\ref{thm-calabi-clifford-cpn} is postponed till
Section~\ref{sect-pf-thm-calabi-clifford-cpn}.

\section{Proof of Theorem~\ref{thm-continuum-quasim-ball-1}}
\label{subsection-pf-thm-continuum-quasim-ball} Let
$B^{2n}(\frac{1}{\sqrt{\pi}}) \subset \C^n$ be the closed
symplectic ball of radius $\frac{1}{\sqrt{\pi}}$. Define an
embedding $\vartheta_\delta :
\textnormal{Int\,}B^{2n}(\frac{1}{\sqrt{\pi}}) \to \C P^n$,
$0<\delta\leq 1$, as: $$ \vartheta_{\delta}: (w_1,\ldots, w_n)
\mapsto \Bigl[ \Bigl( \frac{1}{\pi} - \sum_{i=1}^n \delta |w_i|^2
\Bigr)^{\frac{1}{2}} : \sqrt{\delta} w_1 :\ldots : \sqrt{\delta}
w_n \Bigr],$$ where $w_1,\ldots ,w_n$ are the standard complex
coordinates on $\C^n$. Each embedding $\vartheta_\delta$ is
conformally symplectic: it maps the symplectic form $\omega$ on
$\C P^n$ to $\delta \omega_B$, where $\omega_B$ is the standard
symplectic form on the ball. In particular,
$\vartheta:=\vartheta_1$ is a genuine symplectic embedding, its
image is the complement of a projective hyperplane in $\C P^n$.
The ball $\vartheta
(\textnormal{Int\,}B^{2n}(\frac{1}{\sqrt{\pi}})) \subset \C P^n$
is invariant under the $T^n{\hbox{\rm -action}}$ on $\C P^n$ --
the corresponding $T^n{\hbox{\rm -action}}$ on $\textnormal{Int\,}
B^{2n}(\frac{1}{\sqrt{\pi}})$ is described by the formula: $$
(e^{2\pi s_1},\ldots, e^{2\pi s_n}):\ \ (w_1,\ldots, w_n)\mapsto
(e^{2\pi s_1} w_1 ,\ldots , e^{2\pi s_n} w_n). $$ The orbits of
the action lying inside the ball are the tori
\[
T (r_1, \ldots, r_n) := \{ (w_1,\ldots, w_n)\in\C^n \, | \ \ |w_i|^2 =
r_i,\ i=1,\ldots, n\},
\]
where $r_1, \ldots, r_n$ is a sequence of non-negative numbers
satisfying $r_1 +\ldots +r_n\leq 1/\pi$. In fact $\vartheta_\delta$
maps each $T (r_1, \ldots, r_n)\subset
\textnormal{Int\,}B^{2n}(\frac{1}{\sqrt{\pi}})$ into the torus
$\Phi^{-1} (\pi \delta r_1,\ldots, \pi \delta r_n)$.  In particular,
\begin{equation}
   \label{eq-iden}
   \vartheta_\delta\;
   \Biggl( T \Bigl(\frac{1}{\delta\pi
     (n+1)}, \ldots, \frac{1}{\delta\pi
     (n+1)}\Bigr)\Biggr) = \Phi^{-1} (p_{clif})
\end{equation}
as long as $\delta$ is sufficiently close to $1$.

Put $M = B^{2n}(\frac{1}{\sqrt{\pi}})$, $G=\Ham (M)$. The
conformally symplectic embeddings $\vartheta_{\delta}: M \to \C
P^n$, $0<\delta\leq 1$, induce monomorphisms
$\vartheta_{\delta,\ast} : G \to \Ham (\C P^n)$.  Let $\mu:G \to
\mathbb{R}$ be a Calabi quasimorphism on $\Ham (\C P^n)$ which is
continuous with respect to the Hofer metric (see
Theorem~\ref{thm-existence}).  Then each $\mu_{\delta}
=\delta^{-n-1}\cdot \mu \circ \vartheta_{\delta,\ast}$ is a Calabi
quasimorphism on $G$.

Fix $\delta_0 < 1$ sufficiently close to 1. We will show now that
every finite collection of quasimorphisms of the form $\mu_{\delta}$,
where $\delta \in [\delta_0, 1]$, is linearly independent over $\R$.
Indeed, take any compactly supported Hamiltonian of the form $F =
\tilde{F} (\pi |w|^2)$ on $M$, where $\tilde{F}: \R\to\R$ is a smooth
function supported inside $[-1,1]$ and $w= (w_1,\ldots, w_n)$ are the
complex coordinates as above.
It follows from Theorem~\ref{thm-calabi-clifford-cpn} combined with
formula (\ref{eq-iden}) that for any $\delta \in [\delta_0,1]$
\begin{equation}
   \label{eqn-quasims-evaluated} \mu_{\delta} (\phi_F) =
   \int_{B^{2n}(\frac{1}{\sqrt{\pi}})} F\omega_B^n -
   \delta^{-n-1} \tilde{F} {\Big(}\frac{\delta^{-1}n}{n+1}{\Big )},
\end{equation}
where $\omega_B$ is the symplectic form on
$B^{2n}(\frac{1}{\sqrt{\pi}})$. This immediately yields the linear
independence and proves the theorem for $M =
B^{2n}(\frac{1}{\sqrt{\pi}})$. Finally observe that in the standard
symplectic $\R^{2n}$ balls of different radii are conformally
symplectomorphic and thus their groups of Hamiltonian diffeomorphisms
are isomorphic. This completes the proof.  \Qed

\section{Proof of  Theorem~\ref{thm-calabi-clifford-cpn}}

\label{sect-pf-thm-calabi-clifford-cpn}

We start with the following elementary lemma.

\begin{lemma}
   \label{lem-non-clifford-torus-displaceable} If $p\in\Delta$,
   $p\neq p_{clif}$, then the subset $\Phi^{-1} (p) \subset \C P^n$
   can be displaced from itself by a Hamiltonian diffeomorphism of $\C
   P^n$.
\end{lemma}

Together with Theorem~\ref{thm-cliff-torus-not-displaceable},
which will be proven below, the lemma shows that the Clifford
torus is the only non-displaceable torus among all the Lagrangian
and isotropic tori $\Phi^{-1} (p) \subset \C P^n$, $p\in\Delta$.

\bigskip
\noindent {\bf Proof of the lemma.}  The action of $T^n$ on $\C
P^n$ extends to an action of $U (n+1)\supset T^n$ on $\C P^n$ --
the latter is induced by the natural action of $U (n+1)$ on
$\C^{n+1}$. A permutation of homogeneous coordinates on $\C P^n$
can be viewed as the action on $\C P^n$ of some element from $U
(n+1)$ and thus it represents a Hamiltonian diffeomorphism of $\C
P^n$. On the other hand, such permutations map each orbit of the
$T^n{\hbox{\rm -action}}$ on $\C P^n$ to an orbit of the same
action. Therefore each permutation of homogeneous coordinates on
$\C P^n$ induces an affine transformation of the simplex $\Delta$
which permutes the corresponding vertices. The barycenter
$p_{clif} \in \Delta$ is the only fixed point of the group of
these permutations, hence the lemma is proved.  \Qed

Now we return to the proof of Theorem~\ref{thm-calabi-clifford-cpn}.
Let $\bar{F}$ be a smooth function defined in a neighborhood $V$ of
$\Delta$. Fix $\epsilon > 0$.  Perturb $\bar{F}$ to a function
$\bar{F}^\prime$ so that
\begin{itemize}
  \item{} $\| \bar{F} - \bar{F}^\prime \|_{C^0} < \epsilon$

  \item{} $\bar{F}^\prime$ is constant in a $\gamma$-neighborhood of
   $p_{clif}$, namely $\bar{F}^\prime \equiv \bar{F}^\prime
   (p_{clif})$ in that neighborhood.  (Here $\gamma > 0$ is a
   parameter.)


\end{itemize}

Let $U_0,\ldots, U_m \subset V$ be a finite cover of $\Delta$ with
Euclidean balls of radius $\gamma/2$ such that

\medskip
\noindent (i) $p_{clif}\in U_0$;

\smallskip
\noindent (ii) $p_{clif}\notin {\textrm Closure} (U_j)$ for any
$j\geq 1$.

\medskip Using a partition of unity write
\[
\bar{F}^\prime - \bar{F}^\prime (p_{clif}) = \bar{F}_1 +\ldots +
\bar{F}_m,
\]
where ${\rm supp}\ \bar{F}_j\subset U_j$ for all $j \geq 1$.

Put $F^\prime = \Phi^\ast\bar{F^\prime}$, $F_j = \Phi^\ast \bar{F}_j
$. Note that the Poisson brackets of all these functions vanish since
all fibers of $\Phi$ are Lagrangian or isotropic submanifolds of $\C
P^n$.  Therefore the diffeomorphisms $\phi_{F_j}$ pairwise commute,
hence
$$\phi_{F^\prime} = \phi_{F_1}\circ ... \circ \phi_{F_m }.$$
Since the
restriction of any homogeneous quasimorphism on an abelian subgroup is
a homomorphism we conclude that
\begin{equation}
   \label{sum}
   \mu (\phi_{F^\prime}) =\mu( \phi_{F_1})+ ... +\mu( \phi_{F_m }).
\end{equation}

If $\gamma$ is small enough, the condition (ii) above and
Lemma~\ref{lem-non-clifford-torus-displaceable} guarantee that the
sets $\Phi^{-1} (U_j)$, $j\geq 1$, are displaceable -- they are
neighborhoods of displaceable sets.  Therefore
\[
\mu (\phi_{F_j}) = \int_{\C P^n} F_j\omega^n,\ j\geq 2,
\]
because $\mu$ is a Calabi quasimorphism.  Substituting this into
equation (\ref{sum}) we obtain that
\[
\mu (\phi_{F^\prime}) = \sum_{j=1}^m \int_{\C P^n} F_j\omega^n =
\int_{\C P^n} F^\prime\omega^n - \bar{F}^\prime (p_{clif}).
\]
Since this is true for any $\epsilon > 0$ and $\mu (\phi_F)$ depends
continuously on $F$ (recall that $\mu$ is Lipschitz with respect to
the Hofer metric) this proves Theorem~\ref{thm-calabi-clifford-cpn}.
\Qed

\section{Proof of Theorem~\ref{thm-cliff-torus-not-displaceable}}

\label{pf-thm-cliff-torus-not-displaceable}

\noindent If the Clifford torus is displaceable by a Hamiltonian
diffeomorphism then so is a sufficiently small open neighborhood
$U$ of it. Thus the restriction of the Calabi quasimorphism $\mu$
to the group $G_U$ of Hamiltonian diffeomorphisms supported inside
$U$ is the classical Calabi homomorphism. On the other hand,
Theorem~\ref{thm-calabi-clifford-cpn} states that if $\bar{F}
(p_{clif}) \neq 0$ then the values of $\mu$ and of the Calabi
homomorphism on the Hamiltonian diffeomorphism generated by
$F=\Phi^\ast \bar{F}$ do {\it not} coincide. Hence the Clifford
torus is not displaceable. This proves
Theorem~\ref{thm-cliff-torus-not-displaceable}.  \Qed

\section{An outlook}

\label{sec-outl} It would be interesting to find a more general
setting where one can derive Lagrangian intersection results from
the pure existence of a Calabi quasimorphism which is continuous
with respect to the Hofer metric. Our proof of the fact that the
Clifford torus is non-displaceable (Theorem
\ref{thm-cliff-torus-not-displaceable}) gives an example of this
kind. Below we present some speculations in this direction.

Let $(M^{2n},\omega)$ be a symplectic manifold whose group of
Hamiltonian diffeomorphisms carries such a quasimorphism, say
$\mu$. Assume that for a simplicial complex $\Delta$ there exists
a map $\Phi: M \to \Delta$ with the following property: the
preimage $\Phi^{-1}(p)$ of any $p \in \Delta$ is an isotropic
submanifold of $M$.  Let $ \cal F$ be a class of "smooth"
functions $\bar{F} : \Delta \to \R$ so that the pull-back $F =
\Phi^*\bar{F}$ is smooth.  We assume that $\cal F$ is
"sufficiently ample", in particular there exist partitions of
unity associated to arbitrarily small coverings.

We say that a point $p \in \Delta$ is displaceable if the preimage
$\Phi^{-1}(p)$ can be displaced by a Hamiltonian isotopy and
non-displaceable otherwise.  Arguing exactly as in the proof of
Theorem \ref{thm-calabi-clifford-cpn} we get the following.

\medskip
\noindent
{\bf A.} Assume that there is a point $p_* \in \Delta$ such that the
set $\Delta \setminus \{p_*\}$ consists of displaceable points.  Then
for every smooth function $\bar F \in \cal F$
\begin{equation}
   \label{eq-fant-1}
   \mu (\phi_{F}) = \int_M F \cdot \omega^n - \bar{F}(p_*).
\end{equation}
This immediately yields (cf. the proof of Theorem
\ref{thm-cliff-torus-not-displaceable}) that $p_*$ is
non-displace\-able. As a logical corollary we get that
non-displaceable points do exist!

\medskip
\noindent
{\bf B.} The correspondence
$$\bar{F} \mapsto \mu(\phi_F)$$
is a {\it linear functional} on the
space of smooth functions $\cal F$.  Indeed, all Poisson brackets
$\{F,H\}$ with $F = \Phi^*\bar{F}, H = \Phi^*\bar{H}$ vanish, and
hence the Hamiltonian diffeomorphisms $\phi_F$ and $\phi_H$ commute.
The desired linearity follows from the fact that $\mu$ restricted to
an abelian subgroup is a homomorphism.

\medskip
\noindent
{\bf C.} Denote by $\Sigma \subset \Delta$ the subset of all
non-displaceable points. Assume first that a function $\bar{F} \in
\cal F$ vanishes on $\Sigma$. Using partitions of unity and the Calabi
property of $\mu$ one readily shows that
$$\mu (\phi_{F}) = \int_M F \cdot \omega^n .$$
This in turn yields the
following generalization of formula (\ref{eq-fant-1}): there exists a
measure, say, $\sigma$ on $\Sigma$ so that
\begin{equation}
   \label{eq-fant-2}
   \mu (\phi_{F}) = \int_M F \cdot \omega^n -
   \int_{\Sigma}\bar{F}\cdot d\sigma
\end{equation}
for every $\bar{F} \in \cal{F}$.  For instance in the setting of the
previous sections $\sigma$ is the Dirac measure concentrated at the
value of the moment map corresponding to the Clifford torus.

\medskip

Let us list some specific examples where in our opinion it would be
interesting to test these suggestions. Note that the existence of a
Calabi quasimorphism is meanwhile established for complex projective
spaces, complex Grassmannians and more generally for spherically
monotone symplectic manifolds with semi-simple quantum homology
algebra -- see \cite{EP}.  An extra complication is due to the fact
that in general a Calabi quasimorphism is defined not on the group
$\Ham (M,\omega)$ itself but on its universal cover.

\begin{exam}
   \label{ex-1}
   {\rm Let $\Phi: M \to \Delta$ be the moment map associated to a
     Hamiltonian action of a half-dimensional torus.  Is it true that
     the measure $\sigma$ from formula (\ref{eq-fant-2}) above is
     again the Dirac measure at the barycenter of $\Delta$?

     As an illustration consider the direct product $M = \C
     P^{n_1}\times\ldots\times \C P^{n_l}$ equipped with the monotone
     symplectic structure $\omega = \omega_1 \oplus \ldots \oplus
     \omega_l$.  Monotone means that $[\omega] \in H^2(M)$ is a
     positive multiple of the first Chern class $c_1$ of the (complex)
     tangent bundle of $M$. More explicitly, if we denote by $h_i \in
     H^2({\mathbb{C}}P^{n_i}; \mathbb{Z})$ the positive generator,
     monotonicity means that $[\omega] = \lambda \bigl( (n_1 + 1) h_1
     \oplus \ldots \oplus (n_l+1) h_l \bigr)$ for some $\lambda > 0$.
     Consider the torus action on $M$ defined as the direct product of
     the standard torus actions on the factors.  Denote by $T_{clif}$
     the Lagrangian torus in $M$ which is the product of the Clifford
     tori in $\C P^{n_i}$, $i=1, \ldots, l$. A simple modification of
     our arguments in the previous sections shows that the measure
     $\sigma$ is indeed the Dirac measure concentrated at the value of
     the moment map corresponding to the Clifford torus.  As a
     consequence, $T_{clif}$ cannot be displaced from itself by a
     Hamiltonian isotopy.  }
\end{exam}

\begin{exam}
   \label{ex-2}
   {\rm Let $(M,\omega)$ be a projective algebraic manifold endowed
     with the Fubini-Study form.  It follows from the theory of
     Lagrangian skeletons \cite{Biran} that $M$ carries {\it many}
     singular foliations by Lagrangian tori and thus provides a
     natural playground for our speculation. Note that in this
     situation a careful description of the space of leaves $\Delta$
     is already a non-trivial problem. A potential outcome would be
     the existence of a non-displaceable Lagrangian torus in a large
     class of symplectic manifolds.}
\end{exam}

\begin{exam}
   \label{ex-3}
   {\rm Take the 2-sphere $S^2$ equipped with an area form.  Every
     Morse function, say $\psi$, on $S^2$ defines a tree $T_{\psi}$,
     called the Reeb graph, whose points correspond to connected
     components of the level sets of $\psi$ (see \cite{EP}).  One has
     a natural map $\Psi: S^2 \to T_{\psi}$ whose preimages define a
     singular foliation by circles of $S^2$.  It turns out that the
     corresponding measure $\sigma$ is the Dirac measure concentrated
     at a special point called the median of the tree (see \cite{EP}).

     Given any triple $(M,\Phi,\Delta)$ as in the examples above,
     define its {\it stabilization} by
     $$(S^2 \times M, \Psi \times \Phi, T_{\psi} \times \Delta).$$
     This gives rise to new interesting examples.  }
\end{exam}

\bigskip
\noindent {\bf Acknowledgements.} We thank Slava Kerner and Felix
Schlenk who read a preliminary version of this paper and made a
number of corrections.
\bigskip

\bibliographystyle{alpha}

\bigskip
\noindent
Paul Biran\\
School of Mathematical Sciences \\
Tel Aviv University \\
Tel Aviv 69978, Israel \\
{\it email}: biran@math.tau.ac.il\\

\bigskip
\noindent
Michael Entov\\
Department of Mathematics\\
Technion -- Israel Institute of Technology\\
Haifa 32000, Israel \\
{\it e-mail}: entov@math.technion.ac.il\\

\bigskip
\noindent
Leonid Polterovich \\
School of Mathematical Sciences \\
Tel Aviv University \\
Tel Aviv 69978, Israel \\
{\it email}: polterov@post.tau.ac.il\\

\end{document}